# FROM THE LONG JUMP RANDOM WALK TO THE FRACTIONAL LAPLACIAN

ENRICO VALDINOCI

ABSTRACT. This note illustrates how a simple random walk with possibly long jumps is related to fractional powers of the Laplace operator.
The exposition is elementary and self-contained.

> ''Le discese ardite
> e le risalite
> su nel cielo aperto
> e poi giù il deserto
> e poi ancora in alto
> con un grande salto...''
> (Lucio Battisti, ''Io vorrei... non vorrei... ma se vuoi...'')

The purpose of this note, which is mainly pedagogical, is to show in a simple, concrete example how singular integrals naturally arise as a continuous limit of discrete, long jump random walks, and to recall a simple description of the integral kernels in terms of the Fourier multipliers.

Singular integrals and nonlocal (especially fractional) operators are a classical topic in harmonic analysis and operator theory [Lan72, Ste70] and they are now becoming impressively fashionable because of their connection with many real-world phenomena.

Indeed, nonlocal operators arise in the thin obstacle problem [Caf79], in optimization [DL76], in finance [CT04], in phase transitions [AB98, ABS98, CSM05, SV08b], in stratified materials [SV08a], in anomalous diffusion [MK00], in crystal dislocation [Tol97], in soft thin films [Kur06], in some models of semipermeable membranes and flame propagation [CRS08], in conservation laws [BKW01], in the ultra-relativistic limit of quantum mechanics [FdlL86], in quasi-geostrophic flows [MT96, Cor98], in multiple scattering [DG75, CK98, GK04], in minimal surfaces [CRS07], in materials science [Bat06] and in water waves [Sto57, Zak68, Whi74, CSS92, CG94, NS94, CW95, dlLP96, CSS97, CN00, GG03, HN05, NT08, dlLV08]. See also [Sil05] for further motivation.

From a probabilistic point of view, such nonlocal operators are related to Lévy processes [Ito84, Ber96, BG99, JMW05]. A naive example in Section 1 will show a possible probabilistic interpretation. Then, in Section 2, we will recall an easy recipe for representing the singular integral kernel in terms of the Fourier symbols. The particular (and particularly interesting) case of the fractional Laplacian will be discussed in Section 3.

The exposition is self-contained and no prerequisite is needed: a basic undergraduate math knowledge will suffice.

In the rest of this note, we will discard all the multiplicative normalizing constants. That is, following a convention typical of the lectures on Fourier analysis, we will write $X = Y$ to mean that there is some normalizing constant $C > 0$ such[1] that $X = CY$.

## 1. LONG JUMP RANDOM WALKS AND SINGULAR INTEGRAL KERNELS

Let $\mathcal{K} : \mathbb{R}^n \to [0, +\infty)$ be even, that is $\mathcal{K}(y) = \mathcal{K}(-y)$ for any $y \in \mathbb{R}^n$, and such that

$$\sum_{k \in \mathbb{Z}^n} \mathcal{K}(k) = 1. \qquad (1)$$

Give a small $h > 0$, we consider a random walk on the lattice $h\mathbb{Z}^n$.

---

It is a pleasure to acknowledge the interesting conversations on the subject of this note with Domenico Marinucci. The author has been supported by MIUR Project Variational Methods and Nonlinear Differential Equations and FIRB Project Analysis and Beyond.

[1] We hope that no reader is bothered by the fact that this convention implies, for instance, that $2\pi = 1$. It seems not to be a joke that on February 5, 1897, the House of Representatives of the State of Indiana unanimously passed a bill which would have supported such a new mathematical truth [Ind97].

We suppose that at any unit of time $\tau$ (which may depend on $h$), a particle jumps from any point of $h\mathbb{Z}^n$ to any other point.

The probability for which a particle jumps from the point $hk \in h\mathbb{Z}^n$ to the point $h\tilde{k}$ is taken to be $\mathcal{K}(k - \tilde{k}) = \mathcal{K}(\tilde{k} - k)$.

Note that, differently from the standard random walk, in this process the particle may experience arbitrarily long jumps, though with a small probability.

We call $u(x,t)$ the probability that our particle lies at $x \in h\mathbb{Z}^n$ at time $t \in \tau\mathbb{Z}$.

Of course, $u(x, t + \tau)$ equals the sum of all the probabilities of the possible positions $x + hk$ at time $t$ weighted by the probability of jumping from $x + hk$ to $x$.

That is,
$$u(x, t + \tau) = \sum_{k \in \mathbb{Z}^n} \mathcal{K}(k) u(x + hk, t).$$

Therefore, recalling the normalization in (1),

(2) $$u(x, t + \tau) - u(x, t) = \sum_{k \in \mathbb{Z}^n} \mathcal{K}(k)\big(u(x + hk, t) - u(x, t)\big).$$

Particularly nice asymptotics are obtained in the case in which $\tau = h^\alpha$ and $\mathcal{K}$ is a homogeneous kernel, say, up to normalization factors

(3) $$\mathcal{K}(y) = |y|^{-(n+\alpha)}, \qquad \text{(for } y \neq 0 \text{ and, say, } \mathcal{K}(0) = 0\text{)},$$

with
$$\alpha \in (0, 2).$$

We observe that (1) holds (again, up to normalization) and

(4) $$\frac{\mathcal{K}(k)}{\tau} = h^n \mathcal{K}(hk).$$

Thus, in this case it is convenient to define
$$\psi(y, x, t) = \mathcal{K}(y)\big(u(x + y, t) - u(x, t)\big)$$

and to use (4) to write (2) as

(5) $$\begin{aligned}
\frac{u(x, t + \tau) - u(x, t)}{\tau} &= \sum_{k \in \mathbb{Z}^n} \frac{\mathcal{K}(k)}{\tau}\big(u(x + hk, t) - u(x, t)\big) \\
&= h^n \sum_{k \in \mathbb{Z}^n} \mathcal{K}(hk)\big(u(x + hk, t) - u(x, t)\big) \\
&= h^n \sum_{k \in \mathbb{Z}^n} \psi(hk, x, t).
\end{aligned}$$

Since the latter is just the approximating Riemann sum of
$$\int_{\mathbb{R}^n} \psi(y, x, t)\, dy,$$

by sending $\tau = h^\alpha \to 0^+$ in (5), that is, by taking the continuous limit of the discrete random walk, we obtain
$$\partial_t u(x, t) = \int_{\mathbb{R}^n} \psi(y, x, t)\, dy$$

that is

(6) $$\partial_t u(x, t) = \int_{\mathbb{R}^n} \frac{u(x + y, t) - u(x, t)}{|y|^{n+\alpha}}\, dy.$$

This shows that a simple random walk with possibly long jumps produces, in the limit, a singular integral with a homogeneous kernel.

We remark that the integral

(7) $$\int_{\mathbb{R}^n} \frac{u(x + y) - u(x)}{|y|^{n+\alpha}}\, dy,$$

which appears in (6) has a singularity when $y = 0$.

However, when $\alpha \in (0,2)$ and $u$ is smooth and bounded, such integral is well defined as principal value, that is as

$$\lim_{\epsilon \to 0^+} \int_{\mathbb{R}^n \setminus B_\epsilon} \frac{u(x+y) - u(x)}{|y|^{n+\alpha}} \, dy.$$

Indeed, $|y|^{-(n+\alpha)}$ is integrable at infinity and

$$\int_{B_1} \frac{\nabla u(x) \cdot y}{|y|^{n+\alpha}} = 0$$

as principal value, because the function $y/|y|^{n+\alpha}$ is odd.

Therefore, we may write the singular integral in (7) as principal value near 0 in the form

$$\int_{B_1} \frac{u(x+y) - u(x) - \nabla u(x) \cdot y}{|y|^{n+\alpha}} \, dy$$

and the latter is a convergent integral near 0 because

$$\frac{|u(x+y) - u(x) - \nabla u(x) \cdot y|}{|y|^{n+\alpha}} \leqslant \frac{\|D^2 u\|_{L^\infty} |y|^2}{|y|^{n+\alpha}} = \frac{\|D^2 u\|_{L^\infty}}{|y|^{n-2+\alpha}}$$

which is integrable near 0.

It is also interesting to write the singular integral in (7) as a weighted second order differential quotient. For this, we observe that, substituting $\tilde{y} = -y$, we have that the integral in (7) equals to

$$(8) \qquad \int_{\mathbb{R}^n} \frac{u(x - \tilde{y}) - u(x)}{|\tilde{y}|^{n+\alpha}} \, d\tilde{y}.$$

Therefore, relabeling $\tilde{y}$ as $y$ in (8), we have that

$$(9) \quad \begin{aligned} & 2 \int_{\mathbb{R}^n} \frac{u(x+y) - u(x)}{|y|^{n+\alpha}} \, dy \\ & = \int_{\mathbb{R}^n} \frac{u(x+y) - u(x)}{|y|^{n+\alpha}} \, dy + \int_{\mathbb{R}^n} \frac{u(x-y) - u(x)}{|y|^{n+\alpha}} \, dy \\ & = \int_{\mathbb{R}^n} \frac{u(x+y) + u(x-y) - 2u(x)}{|y|^{n+\alpha}} \, dy. \end{aligned}$$

The equality obtained in (9) shows that the singular integral in (7) may be written, up to a factor 2, as an average of the second incremental quotient $u(x+y) + u(x-y) - 2u(x)$ against the weight $|y|^{n+\alpha}$.

Such a representation is also useful to remove the singularity of the integral at 0, since, for smooth $u$, a second order Taylor expansion gives that

$$\frac{u(x+y) + u(x-y) - 2u(x)}{|y|^{n+\alpha}} \leqslant \frac{\|D^2 u\|_{L^\infty}}{|y|^{n-2+\alpha}}$$

which is integrable near 0.

It is known [Sil05] that the singular integral in (7) is related to the fractional Laplacian $(-\Delta)^{\alpha/2}$. This relation will be outlined here below (see, in particular (14) and (18) below).

It is also interesting to write the displacement of the above random walk at time $n\tau$, for any $n \in \mathbb{N}$. Namely, if $h\epsilon_j \in h\mathbb{Z}^n$ is the jump performed at time $j\tau$ (that is, the "innovation"), the above discussed random walk is made in such a way that the probability that $\epsilon_j$ equals $k$ is $\mathcal{K}(k)$.

The displacement at time $n\tau$ is then the sum of these innovations, that is

$$\sum_{j=1}^n h\epsilon_j.$$

The $\beta$-moment associated to this process is then

$$\sum_{k \in \mathbb{Z}^n} |k|^\beta \mathcal{K}(k) = \sum_{k \in \mathbb{Z}^n \setminus \{0\}} |k|^{\beta - n - \alpha},$$

which is finite if and only if $\beta < \alpha$.

In the probability theory framework, this is interpreted as the innovation being in the domain of attraction of an "$\alpha$-stable random variable" [ST94].

In particular, the associated variance is not finite, thus reflecting that the process is not Gaussian.

## 2. Kernels and Fourier symbols

Given a "nice" (say, smooth and with fast decay, for simplicity) function $u$, the long jump random walk of Section 1 has lead us to the study of integrals of the type

$$\text{(10)} \quad \int_{\mathbb{R}^n} \big(u(x+y) + u(x-y) - 2u(x)\big)\mathcal{K}(y)\,dy,$$

due to (9).

If we call $\mathcal{L}u$ the integral in (10), one may consider $\mathcal{L}$ a linear operator and look for its "symbol" (or "multiplier") in Fourier space.

That is, if $\mathcal{F}$ denotes the Fourier transform, one may think to write

$$\text{(11)} \quad \mathcal{L}u(x) = \mathcal{F}^{-1}\big(\mathcal{S}\,(\mathcal{F}u)\big),$$

for some function $\mathcal{S} : \mathbb{R}^n \to \mathbb{R}$.

The interesting fact is that $\mathcal{K}$ and $\mathcal{S}$ are related as follows:

$$\text{(12)} \quad \mathcal{S}(\xi) = \int_{\mathbb{R}^n} \big(\cos(\xi \cdot y) - 1\big)\mathcal{K}(y)\,dy,$$

up to normalization factors.

To check that (12) holds, one simply Fourier transforms (11) in the variable $x$, calling $\xi$ the corresponding frequency variable: making use of (10) one obtains

$$\begin{aligned}
\mathcal{S}(\xi)\,(\mathcal{F}u)(\xi) &= \mathcal{F}(\mathcal{L}u) \\
&= \mathcal{F}\left(\int_{\mathbb{R}^n} \big(u(x+y) + u(x-y) - 2u(x)\big)\mathcal{K}(y)\,dy\right) \\
&= \int_{\mathbb{R}^n} \Big(\mathcal{F}\big(u(x+y) + u(x-y) - 2u(x)\big)\Big)\mathcal{K}(y)\,dy \\
&= \int_{\mathbb{R}^n} (e^{i\xi \cdot y} + e^{-i\xi \cdot y} - 2)\,(\mathcal{F}u)(\xi)\,\mathcal{K}(y)\,dy \\
&= \int_{\mathbb{R}^n} (e^{i\xi \cdot y} + e^{-i\xi \cdot y} - 2)\,\mathcal{K}(y)\,dy\,(\mathcal{F}u)(\xi) \\
&= 2\int_{\mathbb{R}^n} \big(\cos(\xi \cdot y) - 1\big)\,\mathcal{K}(y)\,dy\,(\mathcal{F}u)(\xi),
\end{aligned}$$

proving (12).

## 3. The fractional Laplacian

The fractional Laplacian may be naturally introduced in the Fourier space.

Indeed, one has that

$$\partial_j u = \mathcal{F}^{-1}\big(i\xi_j(\mathcal{F}u)\big)$$

and therefore

$$-\Delta u = \mathcal{F}^{-1}\big(|\xi|^2(\mathcal{F}u)\big).$$

Thus, it is natural to define, for $\alpha \in (0, 2)$,

$$\text{(13)} \quad (-\Delta)^{\alpha/2}u = \mathcal{F}^{-1}\big(|\xi|^\alpha(\mathcal{F}u)\big).$$

It is known [Lan72, Ste70] that such a fractional Laplacian may be also represented as the principal value of singular integral, namely

$$\text{(14)} \quad (-\Delta)^{\alpha/2}u = \int_{\mathbb{R}^n} \frac{u(x) - u(y)}{|x - y|^{n+\alpha}}\,dy,$$

up to normalizing constants – again, the above integral is intended in the principal value sense.

Notice that, by (9), one can also write (14) as

$$(-\Delta)^{\alpha/2}u = -\int_{\mathbb{R}^n} \frac{u(x+y) + u(x-y) - 2u(x)}{|y|^{n+\alpha}}\,dy,$$

up to normalizing factors.

We give here a simple proof of the equivalence between the definitions in (13) and in (14).

For this, we observe that, in the notation of (10) and (11), we may write (13) and (14) as
$$\mathcal{S}(\xi) = |\xi|^\alpha \text{ and } \mathcal{K}(y) = -|y|^{-(n+\alpha)}.$$
Therefore, by (12), such equivalence boils down to prove that

(15)
$$|\xi|^\alpha = \int_{\mathbb{R}^n} \frac{1 - \cos(\xi \cdot y)}{|y|^{n+\alpha}} \, dy.$$

To prove (15), first observe that, if $\zeta = (\zeta_1, \ldots, \zeta_n) \in \mathbb{R}^n$, we have
$$\frac{1 - \cos \zeta_1}{|\zeta|^{n+\alpha}} \leqslant \frac{|\zeta_1|^2}{|\zeta|^{n+\alpha}} \leqslant \frac{1}{|\zeta|^{n-2+\alpha}}$$
near $\zeta = 0$, therefore

(16) $\int_{\mathbb{R}^n} \frac{1 - \cos \zeta_1}{|\zeta|^{n+\alpha}} \, d\zeta$ is finite and positive.

We now consider the function
$$\mathcal{J}(\xi) = \int_{\mathbb{R}^n} \frac{1 - \cos(\xi \cdot y)}{|y|^{n+\alpha}} \, dy.$$

We have that $\mathcal{J}$ is rotationally invariant, that is

(17) $$\mathcal{J}(\xi) = \mathcal{J}(|\xi|e_1).$$

Indeed, if $n = 1$, then one easily checks that $\mathcal{J}(-\xi) = \mathcal{J}(\xi)$, proving (17) in this case.

When $n \geqslant 2$, we consider a rotation $R$ for which
$$R(|\xi|e_1) = \xi$$
and we denote by $R^T$ its transpose. We obtain, via the substitution $\tilde{y} = R^T y$,

$$\begin{aligned}
\mathcal{J}(\xi) &= \int_{\mathbb{R}^n} \frac{1 - \cos\left((R(|\xi|e_1)) \cdot y\right)}{|y|^{n+\alpha}} \, dy \\
&= \int_{\mathbb{R}^n} \frac{1 - \cos\left((|\xi|e_1) \cdot (R^T y)\right)}{|y|^{n+\alpha}} \, dy \\
&= \int_{\mathbb{R}^n} \frac{1 - \cos\left((|\xi|e_1) \cdot \tilde{y}\right)}{|\tilde{y}|^{n+\alpha}} \, dy \\
&= \mathcal{J}(|\xi|e_1),
\end{aligned}$$

which proves (17).

As a consequence of (17) and (16), the substitution $\zeta = |\xi|y$ gives that

$$\begin{aligned}
\mathcal{J}(\xi) &= \mathcal{J}(|\xi|e_1) \\
&= \int_{\mathbb{R}^n} \frac{1 - \cos(|\xi|y_1)}{|y|^{n+\alpha}} \, dy \\
&= \frac{1}{|\xi|^n} \int_{\mathbb{R}^n} \frac{1 - \cos \zeta_1}{|\zeta/|\xi||^{n+\alpha}} \, d\zeta = |\xi|^\alpha,
\end{aligned}$$

up to normalization factors, hence (15) is proved, thus so is the equivalence between (13) and (14).

We remark that, from (14), the probability density of the limit long jump random walk in (6) may be written as

(18) $$\partial_t u = -(-\Delta)^{\alpha/2} u.$$

## References

[AB98]  Giovanni Alberti and Giovanni Bellettini. A nonlocal anisotropic model for phase transitions. I. The optimal profile problem. *Math. Ann.*, 310(3):527–560, 1998.

[ABS98] Giovanni Alberti, Guy Bouchitté, and Pierre Seppecher. Phase transition with the line-tension effect. *Arch. Rational Mech. Anal.*, 144(1):1–46, 1998.

[Bat06]  Peter W. Bates. On some nonlocal evolution equations arising in materials science. In *Nonlinear dynamics and evolution equations*, volume 48 of *Fields Inst. Commun.*, pages 13–52. Amer. Math. Soc., Providence, RI, 2006.

[Ber96]  Jean Bertoin. *Lévy processes*, volume 121 of *Cambridge Tracts in Mathematics*. Cambridge University Press, Cambridge, 1996.

[BG99]     Tomasz Bojdecki and Luis G. Gorostiza. Fractional Brownian motion via fractional Laplacian. *Statist. Probab. Lett.*, 44(1):107–108, 1999.
[BKW01]    Piotr Biler, Grzegorz Karch, and Wojbor A. Woyczyński. Critical nonlinearity exponent and self-similar asymptotics for Lévy conservation laws. *Ann. Inst. H. Poincaré Anal. Non Linéaire*, 18(5):613–637, 2001.
[Caf79]    L. A. Caffarelli. Further regularity for the Signorini problem. *Comm. Partial Differential Equations*, 4(9):1067–1075, 1979.
[CG94]     Walter Craig and Mark D. Groves. Hamiltonian long-wave approximations to the water-wave problem. *Wave Motion*, 19(4):367–389, 1994.
[CK98]     David Colton and Rainer Kress. *Inverse acoustic and electromagnetic scattering theory*, volume 93 of *Applied Mathematical Sciences*. Springer-Verlag, Berlin, second edition, 1998.
[CN00]     Walter Craig and David P. Nicholls. Travelling two and three dimensional capillary gravity water waves. *SIAM J. Math. Anal.*, 32(2):323–359 (electronic), 2000.
[Cor98]    Diego Cordoba. Nonexistence of simple hyperbolic blow-up for the quasi-geostrophic equation. *Ann. of Math. (2)*, 148(3):1135–1152, 1998.
[CRS07]    Luis Caffarelli, Jean-Michel Roquejoffre, and Ovidiu Savin. Non-local minimal surfaces. *In preparation*, 2007.
[CRS08]    Luis Caffarelli, Jean-Michel Roquejoffre, and Yannick Sire. Free boundaries with fractional Laplacians. *In preparation*, 2008.
[CSM05]    Xavier Cabré and Joan Solà-Morales. Layer solutions in a half-space for boundary reactions. *Comm. Pure Appl. Math.*, 58(12):1678–1732, 2005.
[CSS92]    W. Craig, C. Sulem, and P.-L. Sulem. Nonlinear modulation of gravity waves: a rigorous approach. *Nonlinearity*, 5(2):497–522, 1992.
[CSS97]    Walter Craig, Ulrich Schanz, and Catherine Sulem. The modulational regime of three-dimensional water waves and the Davey-Stewartson system. *Ann. Inst. H. Poincaré Anal. Non Linéaire*, 14(5):615–667, 1997.
[CT04]     Rama Cont and Peter Tankov. *Financial modelling with jump processes*. Chapman & Hall/CRC Financial Mathematics Series. Chapman & Hall/CRC, Boca Raton, FL, 2004.
[CW95]     Walter Craig and Patrick A. Worfolk. An integrable normal form for water waves in infinite depth. *Phys. D*, 84(3-4):513–531, 1995.
[DG75]     J. J. Duistermaat and V. W. Guillemin. The spectrum of positive elliptic operators and periodic bicharacteristics. *Invent. Math.*, 29(1):39–79, 1975.
[DL76]     G. Duvaut and J.-L. Lions. *Inequalities in mechanics and physics*. Springer-Verlag, Berlin, 1976. Translated from the French by C. W. John, Grundlehren der Mathematischen Wissenschaften, 219.
[dlLP96]   R. de la Llave and P. Panayotaros. Gravity waves on the surface of the sphere. *J. Nonlinear Sci.*, 6(2):147–167, 1996.
[dlLV08]   Rafael de la Llave and Enrico Valdinoci. Symmetry for a Dirichlet-Neumann problem arising in water waves. *Preprint*, 2008.
[FdlL86]   C. Fefferman and R. de la Llave. Relativistic stability of matter. I. *Rev. Mat. Iberoamericana*, 2(1-2):119–213, 1986.
[GG03]     Günter K. Gächter and Marcus J. Grote. Dirichlet-to-Neumann map for three-dimensional elastic waves. *Wave Motion*, 37(3):293–311, 2003.
[GK04]     Marcus J. Grote and Christoph Kirsch. Dirichlet-to-Neumann boundary conditions for multiple scattering problems. *J. Comput. Phys.*, 201(2):630–650, 2004.
[HN05]     Bei Hu and David P. Nicholls. Analyticity of Dirichlet-Neumann operators on Hölder and Lipschitz domains. *SIAM J. Math. Anal.*, 37(1):302–320 (electronic), 2005.
[Ind97]    The State of Indiana. The Indiana Pi Bill, 1897. http://www.agecon.purdue.edu/crd/Localgov/Secondbill.htm.
[Ito84]    K. Ito. *Lectures on stochastic processes*, volume 24 of *Tata Institute of Fundamental Research Lectures on Mathematics and Physics*. Distributed for the Tata Institute of Fundamental Research, Bombay, second edition, 1984. Notes by K. Muralidhara Rao.
[JMW05]    Benjamin Jourdain, Sylvie Méléard, and Wojbor A. Woyczynski. A probabilistic approach for nonlinear equations involving the fractional Laplacian and a singular operator. *Potential Anal.*, 23(1):55–81, 2005.
[Kur06]    Matthias Kurzke. A nonlocal singular perturbation problem with periodic well potential. *ESAIM Control Optim. Calc. Var.*, 12(1):52–63 (electronic), 2006.
[Lan72]    N. S. Landkof. *Foundations of modern potential theory*. Springer-Verlag, New York, 1972. Translated from the Russian by A. P. Doohovskoy, Die Grundlehren der mathematischen Wissenschaften, Band 180.
[MK00]     Ralf Metzler and Joseph Klafter. The random walk's guide to anomalous diffusion: a fractional dynamics approach. *Phys. Rep.*, 339(1):77, 2000.
[MT96]     Andrew J. Majda and Esteban G. Tabak. A two-dimensional model for quasigeostrophic flow: comparison with the two-dimensional Euler flow. *Phys. D*, 98(2-4):515–522, 1996. Nonlinear phenomena in ocean dynamics (Los Alamos, NM, 1995).
[NS94]     P. I. Naumkin and I. A. Shishmarëv. *Nonlinear nonlocal equations in the theory of waves*. American Mathematical Society, Providence, RI, 1994. Translated from the Russian manuscript by Boris Gommerstadt.
[NT08]     David P. Nicholls and Mark Taber. Joint analyticity and analytic continuation of Dirichlet-Neumann operators on doubly perturbed domains. *J. Math. Fluid Mech.*, 10(2):238–271, 2008.
[Sil05]    Luis Silvestre. *Regularity of the obstacle problem for a fractional power of the Laplace operator*. PhD thesis, University of Texas at Austin, 2005.
[ST94]     Gennady Samorodnitsky and Murad S. Taqqu. *Stable non-Gaussian random processes*. Stochastic Modeling. Chapman & Hall, New York, 1994. Stochastic models with infinite variance.

[Ste70] Elias M. Stein. *Singular integrals and differentiability properties of functions*. Princeton University Press, Princeton, N.J., 1970. Princeton Mathematical Series, No. 30.
[Sto57] J. J. Stoker. *Water waves: The mathematical theory with applications*. Pure and Applied Mathematics, Vol. IV. Interscience Publishers, Inc., New York, 1957.
[SV08a] Ovidiu Savin and Enrico Valdinoci. Elliptic PDEs with fibered nonlinearities. *To appear in J. Geom. Anal.*, 2008.
[SV08b] Yannick Sire and Enrico Valdinoci. Fractional Laplacian phase transitions and boundary reactions: a geometric inequality and a symmetry result. *Preprint*, 2008.
[Tol97] J. F. Toland. The Peierls-Nabarro and Benjamin-Ono equations. *J. Funct. Anal.*, 145(1):136–150, 1997.
[Whi74] G. B. Whitham. *Linear and nonlinear waves*. Wiley-Interscience [John Wiley & Sons], New York, 1974. Pure and Applied Mathematics.
[Zak68] V.E. Zakharov. Stability of periodic waves of finite amplitude on the surface of a deep fluid. *Zh. Prikl. Mekh. Tekh. Fiz.*, 9:86–94, 1968.

Dipartimento di Matematica, Università di Roma Tor Vergata, Via della Ricerca Scientifica, 1, I-00133 Roma (Italy)

*E-mail address*: enrico@math.utexas.edu